\documentstyle[12pt,twoside]{article}

\makeatletter
\renewcommand{\@oddhead}{reducing to Krein-isometry \hfill \thepage}
\renewcommand{\@evenhead}{\thepage \hfill Sergej A. Choro\v savin }
\renewcommand{\@oddfoot}{}
\renewcommand{\@evenfoot}{}
\makeatother

{\vspace*{2ex}\par {#1}\quad\it }{\medskip\par }

\newenvironment{Thm}[2]{\par\addvspace{\bigskipamount}{\bf #1#2}\it }%
{\par\addvspace{\bigskipamount} }

\newenvironment{Theorem}[1]{\begin{Thm}{Theorem}{#1}}{\end{Thm} }

\newenvironment{Remark}[1]{\begin{Thm}{Remark}{#1}}{\end{Thm} }

\newenvironment{Proof}{\par\addvspace{\bigskipamount} {\sc Proof}}%
{\par\hspace*{\fill}$\Box$ \par\addvspace{\smallskipamount} }

\author{S.A.~Choro\v{s}avin}  
\title{A Reducing of the Invariant Semidefinite Subspace Problem
 for Kre\u{\i}n Noncontraction to such a Problem
for Kre\u{\i}n Isometry }
\date{}

\begin{document}
\maketitle

\begin{abstract}
  Theorem.
If every J-isometry has nontrivial positive invariant subspace
then every J-noncontraction has such a subspace.

  Theorem.
If every J-binoncontractive J-isometry has maximal positive invariant subspace
then every J-noncontraction has such a subspace.

\end{abstract}

{\bf 0.}

   Let $H$ be a Hilbert space equipped with the fixed orthogonal decomposition
$ H = H_-  + H_+ $. The indefinite inner product
           $$x,y \in H \mapsto <x,y> \in C$$
is introduced now by the formula
          $$<x,y>:=(x,Jy);\qquad  J:=P_-+P_+ ;$$
where $( , )$ is the symbol of the usual Hilbert scalar product and
$P_{\pm} $ are
the orthoprojectors of $H$ onto $H_{\pm }$. The pair $H,< , >$
is said to be {\it Krein space} or {\it J-space}\/.

     A linear bounded $T: H \to H$ is said to be {\it J-noncontraction}
or {\it $< , >$ -noncontraction} 
or {\it Krein noncontraction\/}, 
iff
              $$<Tx,Tx> \quad \geq \quad <x,x> \qquad   (x\in H) $$
Particular cases of Krein noncontraction are {\it J-isometry ($<,>$-isometry,
Krein-isometry)\/}:
$$
                  <Vx,Vx>=<x,x>   \qquad (x\in H)
$$
{\it Krein-unitary\/}:
$$
      <Ux,Ux>\;=\;<x,x>,\qquad <U^\dagger x,U^\dagger x>\;=\;<x,x>
 \qquad (x\in H)
$$
and {\it Krein binoncontraction\/}:
$$
     <Tx,Tx>\;\geq\;<x,x>,\qquad <T^\dagger x,T^\dagger x>\;\geq \;<x,x>
                                                \qquad (x\in H)
$$
hereinafter $T^\dagger$ is $<,>$-adjoint of $T$ (Note,
$T^\dagger = JT^* J$, where $T^*$ is the Hilbert adjoint of $T$).

The traditional question of the theory of linear operators is
the question of the existense of a (nontrivial) invariant subspace.
In the case of indefinite inner product spaces it is of interest
a special kind of subspaces:

A subspace $L\subset H$ is said to be {\it positive (negative)\/} 
iff
$$
     <x,x>\; \geq\; 0 \quad (\leq 0)\quad \mbox{for any} \; x\in L.
$$
     A subspace $L\;\subset H$ is said to be 
{\it maximal positive (negative)\/}
iff it is positive (negative) and maximal through such of subspaces in the
sence of the set theory.

    {\it In order that the subspace $L\subset H$ should be maximal positive,
a necessary and sufficient condition is that there should exist a linear
$K:H_+\to H_-$ such that }
$$
\|K\|\leq 1 \mbox{ and }
    L=\{x_+ +Kx_+ | x_+\in H_+\}.$$ 
Such an operator $K$ is said to be
{\it angular operator\/} 
of $L$ and it is unique (by fixed $L$). The structure
of cloused positive subspace $L$ is analogous:$L=\{x_++Kx_+|x_+\in L_+\}$,
where $L_+$ is a cloused linear subspace of $H_+$ and $K:L_+\to H_-$ is
linear with $\|K\|\leq 1$.

     We are interested in the existence  of maximal positive invariant
subspace for Krein noncontraction.

    At first sight the question for noncontraction seems to be more general
than the one for isometry and the last question seems to be more general than
the question for Krein unitary operator.
    We shall show that it is not entirely the case.

{\bf 1.}

\begin{Theorem}{ 1.} 

If every  Krein  isometry has a nontrivial
positive invariant subspace \,(resp. maximal positive  invariant subspace),
\, then every Krein  noncontraction has such a subspace.
\end{Theorem}
\begin{Proof}. 
Let us consider a denumerable Hilbert direct sume
$$
      \hat H := H\oplus H\oplus ...=\bigoplus _{n=1}^\infty H
$$
with the canonical projection $p:=\hat H \to H$
$$
     p:x_1\oplus x_2\oplus ... \mapsto x_1
$$
canonical embedding $j:H\to \hat H$
$$
      j:x\mapsto x\oplus 0\oplus 0 ...
$$
and introduce operators $\hat P_\pm :\hat H\to \hat H,\quad \hat J:
\hat H\to \hat H  $   by the formulas:
$$
   \hat P_+:=P_+\oplus 0\oplus 0\oplus ...
$$
$$
   \hat P_-:=P_-\oplus I\oplus I\oplus...
$$
$$
   \hat J \quad :=J\oplus - I\oplus -I\oplus ...
$$
Note: the operators $\hat P_\pm $ are orthoprojections, $\hat P_++\hat P_-=I
_{\hat H}, \hat J =\hat P_++\hat P_-$, and $\hat H$ is a Krein space
with respect to the decomposition
$$
    \hat H=\hat H_++\hat H_-, \quad \hat H_\pm := \hat P_\pm\hat H
$$
    Let $\ll ,\gg$ will denote the correspondent indefinite inner product on
$\hat H$ . We have to remark the next properties of it:
$$
1)\mbox{ if } \ll x,x \gg \ge 0 \mbox{ then } \ll px,px\gg \ge 0; (x\in \hat H)
$$
$$
2)\mbox{ if } \ll x,x\gg \ge 0 \mbox{ and } x\ne 0 \mbox{ then } px\ne 0;
(x\in \hat H)
$$
and hence if $L$ is a nontrivial subspace of $\hat H$, then the $pL$ is a
 nontrivial positive subspace of $H$. Moreover if $L$ is a maximal
positive subspace of $\hat H$ then the $pL$ is a maximal positive subspace
of $H$. To show this it is sufficient to observe that $p\hat P_\pm =P_\pm ,
jP_+=\hat P_+$. Hence if $L$ is a maximal positive with the angular operator
$K:\hat H_+ \to \hat H_-$:
$$
    L =\{ x_+ +Kx_+|x_+\in \hat H_+ \} 
$$
then
$$
   pL=\{ px_++pKx_+|x_+\in \hat H\} 
$$   
$$
    =\{ \tilde x_++pKj\tilde x_+|\tilde x_+ \in H_+ \} ,
$$
$\| pKj\| \leq \| K\| \leq 1 $ and $pL$ is maximal positive with the angular
operator $pKj$.

    Now the rest is fast evident: 

    If $T$ is a Krein noncontraction, then $T^\ast JT-J\geq 0$ and there   
exists the Hilbert square root $ D=(T^\ast JT -J)^{1/2} $. Constructing the
operator $V:\hat H \to \hat H$ as follows
$$
    V:x_1 \oplus x_2 \oplus x_3 \oplus ...\mapsto 
      Tx_1 \oplus Dx_1 \oplus x_2 \oplus x_3 \oplus ...
$$
we obtain
$$
     Tp=pV,\qquad \ll Vx,Vx\gg =\ll x,x\gg \qquad (x\in \hat H)
$$

    Hence $V$ is a Krein isometry and if a subspace $L$ is invariant for
$V$ then the $pL$ is invariant for T.$\Box $
\end{Proof} 

\begin{Remark}. 
If $dimH_\pm =\infty $ both, then $\hat H$ is isomorph
to $H$ as a Krein space. Hence for this case the theorem 1 can be reformulated
so : 
 If every  Krein  isometry $V:H\to H$ has a nontrivial
positive invariant subspace (resp. maximal positive  invariant subspace),
then every Krein noncontraction $T:H\to H$ has such a subspace.$\Box $ 
\end{Remark}

{\bf 2.}
\begin{Theorem}{ 2.} 
  
If every Krein binoncontractive Krein isometry
has a maximal positive invariant subspace, then every Krein isometry has
such a subspace.
\end{Theorem}
 
\begin{Proof}. 
Let $\hat H:=H\oplus H$ be the usual Hilbert
direct sume with canonical projection $p:\hat H \to H$ defined by the formula
$$
     p:x_1 \oplus x_2 \mapsto x_1   \qquad  (x_1,x_2 \in H)
$$
and let $\hat J,\hat V: \hat H \to \hat H $ be linear operators defined as
follows:   $\hat J:=J\oplus I    $  ,

$$
\hat V:= \left(
\begin{array}{cc} 
V & A \\  0 & B^*
\end{array}
\right)
$$
where $V:H\to H$ is a Krein isometry  and $A,B: H\to H$ are linear and 
bounded. The space $\hat H$ is the Krein space with respect to $\hat J$
and if $L$ is a maximal positive subspace of $\hat H$ then
 $L':=p(L\cap H\oplus \{0\}) $ is a maximal positive subspace of the
   space  $H$. If $L$ is invariant for $\hat V$, then $L'$ is invariant
 for $V$.

    We shall demonstrate that there exist $A$ and $B$ such that the operator
$\hat V$ is a Krein binoncontractive Krein isometry.

    Previously we have to make some observations. \\  
1) $J \, -\, VJV^*$ is selajoint and
$$
    Ran\, (J-VJV^*)\; =\; Ran\, (I-VV^+)\; =\; Ker\, V^+=Ker\, V^*J;
$$
let $\, J-VJV^*\, =\, \int \limits_{-\infty}^{+\infty}\lambda \, dE_{\lambda}$
be the spectral decomposition and put
$$
p_+:=\int_{+0}^{\infty}  dE_{\lambda},\hfill 
$$
$$
A:= (J-VJV^*)_+^{1/2}\; =\; \int_{+0}^{\infty}\lambda^{1/2}\, dE_{\lambda} 
\hfill
$$
So $A$ and $p_+$ are selfadjoint, positive and $p_+$ is an orthoprojector onto
$\overline{Ran\, A}$. Hence
$$
Ran\, p_+ = \overline{Ran\, (J-VJV^*)p_+}\subset Ran\, (J-VJV^*)\, =\, 
Ran\, (I -VV^{\dagger})\, 
$$
$$
\qquad \qquad =\, Ker\, V^{\dagger}\, = \, Ker\, V^*J
$$
2) We have $(JV)^*=V^*J$ and hence
$$
H\, =\, Ker\, V^*J+Ran\, JV, \quad Ker\, V^*J \perp Ran\, JV
$$

   Let $JV=u| JV|$ be the correspondent polar decomposition
where $u$ is the partial isometry with final space $Ran\, JV$
and initial space $Ran\, V^*J =Ran\, V^+ =H$; hence, $u$ is an isometry.
But
$$
Ran\, (I-p_+)=Ker\, p_+\, \supset \, (Ker\, V^*J)^{\perp}\, 
=\, Ran\, JV
$$
hence there exist an isometry of the space $H$  with final space 
$Ran\, (I-p_+)$; we define the operator $B$ as such an isometry.

Now we obtain by inmediate computations:
$$
\hat V^* \hat J \hat V = \left(
\begin{array}{cc}
V^*JV & V^*JA \\ AJV & AJA+BB^*
\end{array} \right)
= \left(
\begin{array}{cc}
J & V^*JA \\ AJV & AJA+(I-p_+)
\end{array}
\right)
$$
$$
\hat V \hat J \hat V^* = \left(
\begin{array}{cc}
VJV^* + A^2 & AB\\ B^*A^* & B^*B
\end{array} \right) = \left(
\begin{array}{cc}
VJV^* +A^2 & AB\\ B^*A^* & I 
\end{array}
\right)
$$
Using the identities $A=Ap_+=p_+A=p_+Ap_+$, $B=(I-p_+)B$ we obtain
$V^*JA=V^*Jp_+A=0$ (and hence $AJV=0$), $AB=Ap_+(I-p_+)B$
(and hence $B^*A^*=0$),
$$
\begin{array}{cc}
A^2JA^2&=|J-VJV^*|p_+ Jp_+ |J-VJV^*| \hfill \\
       &=(J-VJV^*)p_+ Jp_+ (J-VJV^*) \hfill \\
       &=p_+ (J-VJV^*)J(J-VJV^*)p_+ \hfill \\
       &=p_+(J-VJV^*)p_+=p_+|J-VJV^*|p_+\hfill \cr
       &=(|J-VJV^*|^{1/2}p_+)p_+ (|J-VJV^*|^{1/2}p_+) \hfill \cr
       &=Ap_+A \hfill
\end{array}
$$
But $A$ invertible on $RanA$, hence $AJA=p_+$ and $ \hat V^*\hat J \hat V\;=\;
p_+$ . 

    Finally we have consequently
$$
A^2=|J-VJV^*|p_+\geq J-VJV^*,\quad VJV^*+A^2\geq J, \quad 
\hat V\hat J\hat V^*\geq \hat J
$$
and hence $\hat V$ is a$ \hat J$-binoncontractive $\hat J$-isometry.
    
\end{Proof} 

\medskip

  The article text is the complete text of the author's report on 15-th
Voronezh Winter Mathematical School  [4].
But in that time the presented construtions and theorems seemed to be rather
curious observations. Now the situattion is changing (see e.g. 
LANL E-print math.DS/9908169
or [5], [6])

\bibliographystyle{unsrt}

\end{document}